 \documentclass{amsart}
\newtheorem{theorem}{Theorem}[section]
\newtheorem{lemma}[theorem]{Lemma}

\newtheorem{definition}[theorem]{Definition}
\newtheorem{example}[theorem]{Example}

\newtheorem{corollary}[theorem]{Corollary}
\newtheorem{proposition}[theorem]{Proposition}

\newtheorem{remark}[theorem]{Remark}

\numberwithin{equation}{section}

\renewcommand{\oddsidemargin}{0.5cm}

\def\C{\mathbb C}
\def\R{\mathbb R}
\def\X{\mathbb X}
\def\Z{\mathbb Z}
\def\Y{\mathbb Y}
\def\Z{\mathbb Z}
\def\N{\mathbb N}
\def\cal{\mathcal}
\def\cD{\cal D}
\def\tD{\tilde{{\cal D}}}
\def\F{\cal F}
\def\tf{\tilde{f}}
\def\tg{\tilde{g}}

\begin{document}
\title[Spectral Theory of Functions \& Applications]{A Spectral Theory of Non-Uniformly Continuous Functions and
the Loomis-Arendt-Batty-Vu Theory on the Asymptotic Behavior of Solutions of Evolution Equations}
\author{Nguyen Van Minh}
\address{Department of Mathematics, University of West Georgia, Carrollton, GA 30118}
\email{vnguyen@westga.edu}

\date{\today}
\begin{abstract}In this paper we present a new
approach to the spectral theory of {\it non-uniformly continuous} functions and a new framework
for the Loomis-Arendt-Batty-Vu theory. Our approach is direct and
free of $C_0$-semigroups, so the obtained results, that extend
previous ones, can be applied to large classes of evolution
equations and their solutions.
\end{abstract}

\maketitle


\section{Introduction}
There is a remarkable theory of the asymptotic behavior of orbits
of bounded one-parameter semigroups of operators on Banach spaces
in the case when the purely imaginary part of the spectrum of the
generator is countable. Although it can be viewed as part of Pure
Functional Analysis and Operator Theory, it has direct
applications to the asymptotic behavior of solutions of evolution
equations. This theory is closely related to Tauberian Theory,
Loomis' theorem in Harmonic Analysis concerning almost periodicity
of functions with countable Beurling spectrum, and a Gelfand's
theorem on the spectrum of $C_0$-groups of isometries.

\medskip
Significant contributions in the theory which we refer to as the
Loomis-Arendt-Batty-Vu theory, were made by L. Loomis, W. Arendt,
C. Batty, Vu Quoc Phong and others. We refer the reader to the
monographs \cite{arebathieneu,nee} and their references for more
systematic information on the theory. The Loomis-Arendt-Batty-Vu
theory is based on early works extending Loomis' theorem to study
almost periodic solutions of evolution equations,
\cite{levzhi,bask,bask2}. In this direction, general results of
the type of Loomis' theorem were proved by applying a Gelfand's
theorem, and the concept of reduced Beurling spectrum (see e.g.
\cite{arebat2,bas,baspry,bask,bask2,levzhi,ruevu}). A far-reaching
ergodic condition, introduced in \cite{ruevu} to study
asymptotically almost periodic functions on the line (see also
\cite{bas}) turns out to be a key tool to present the theory in a
unified framework, especially for asymptotic almost periodicity,
asymptotic stability of solutions of evolution equations on the
half line (see e.g.
\cite{arebat,arebat2,arebathieneu,batneerab2,chifas,delvu,lyuvu,lyuvu2,vu}).

\medskip
There is a "prevalence of the hypothesis of uniform continuity in
the literature", as remarked by Loomis in the Introduction of his
paper \cite{loo}. This remark seems to have been valid until now,
and can be explained by the use of sophisticated tools involving
semigroups of operators. As the Loomis-Arendt-Batty-Vu theory has
direct applications to evolution equations, the requirement that
the equations be well-posed or the solutions be uniformly
continuous for the use of semigroup theory seems to be technical,
and is an obstacle for applications.

\medskip
In this paper we will take an attempt to push forward the
Loomis-Arendt-Batty-Vu theory by introducing a new approach to the
spectral theory of {\it not necessarily uniformly continuous functions}, and a new framework for the main
points of the theory. Our approach is more elementary, free of
$C_0$-semigroups. As a result, the theory with general results can be presented in a way that is easily accessible to readers working in ordinary differential equations with a limited knowledge of the methods of harmonic analysis.

\medskip
We will start the paper with a new approach to the concept of
reduced spectrum of a bounded function on $\R$ or on $\R^+$ that is {\it not necessarily uniformly continuous}. This
approach is based on some simple facts from the ODE that gives
more direct insights into the differentiation operator $\cD :=
d/dt$ on various function spaces. As a result, a new relation
between the reduced spectrum and spectrum of the differentiation
operator $\cD$ is established (Theorem \ref{the spe of tD
Lambda}). We also extend the Gelfand's theorem (Theorem \ref{the gel}) as the key tool to prove Loomis Theorem (Theorem \ref{the
loomis}). Among consequences of the Loomis Theorem we mention
standard ones without the hypothesis of uniform continuity
(Corollaries \ref{cor scalar}, \ref{cor loomis 2}). As an almost
automorphic function may not be uniformly continuous, Corollary
\ref{cor aa} is a typical example of a result which has not been
covered in previous works. We emphasize that a version of Loomis
Theorem with an ergodic condition can be easily derived following
the lines discussed in the next part of the paper. For functions
on the half line, our direct approach to the differentiation
operator $\cD$ on $\R^+$ results in Theorem \ref{the spe dif ope
hl} which is a key tool to prove the main result of the section
(Theorem \ref{the loomis half line}). In the applications of our
results to stability of evolution equations we can free both the
hypotheses of uniform continuity and well-posedness (Corollary
\ref{cor 4.3}, Theorem \ref{the minh}). In particular, our
approach also yields the well known Arendt-Batty-Lyubich-Vu
Theorem (Corollary \ref{the ablv}). In our references we give a
(non-exhausted) list of works concerned with the theory and
related applications of spectral theory of functions to the
asymptotic behavior of evolution equations. (The reader may use
the references of \cite{arebathieneu,nee,vu2} for a more complete
list of references for this paper.)

\medskip
The first version of this paper was announced in the preprint \cite{min} (see also \cite{min2,minngusie} for a related result).
In this new version, we have added some new developments from the literature since
then with appropriate comments as well as corrections and applications. \footnote{After the first version \cite{min} of this paper appeared we learned that in \cite{basgun2} related concepts are discussed. However, the related results stated in \cite{basgun2} seem to be different and incomplete. For example, in addition to the uniform closedness condition which seems to be very strict, Proposition 4.1 and Theorem 4.3 obviously assume that the reduced spectrum of the considered function is non-empty. This condition appears to be significant if the Geldfand Theorem for isometries is not available. }

\bigskip
Before closing this section we would like to list some standard
notations we will use in the paper. Throughout the paper $\R$
denotes the real line, $\R^+$ denotes the half line $[0,\infty )$,
and $\X$ denotes a Banach space over the complex plane $\C$. If
$A$ is a linear operator on a Banach space $\X$, $D(A)$ stands for
its domain; $\sigma(A)$, $\sigma_p(A)$, $\rho (A)$ stand for its
spectrum, point spectrum and resolvent set, respectively. $L(\X)$
stands for the Banach space of all bounded linear operators in
$\X$ with the usual norm $\| \cdot \|$. If $\lambda\in \rho (A)$,
then $R(\lambda ,A)$ denotes the resolvent $(\lambda -A)^{-1}$. In
this paper we will use the following notations:
\begin{enumerate}
\item $c_0$ denotes the Banach space of all numerical sequence
$x=\{ x_n\}_{n=1}^\infty$ such that $\lim_{n\to\infty} x_n =0$,
with sup-norm $\| x\| =\sup_{n\in\N}\| x_n\|$; \item $J$ is either
$\R$ or $\R^+$; \item $BC(J,\X)$ is the space of all $\X$-valued
bounded and continuous functions on $J$; \item $BUC(J,\X)$ is the
space of all $\X$-valued bounded and uniformly continuous
functions on $J$; \item $AP(\X), AP(\R^+,\X)$ are the spaces of
all $\X$-valued almost periodic functions on $\R$, and $\R^+$,
respectively;
 \item $AA(\X), AA(\R^+,\X)$ are the spaces of all
$\X$-valued almost automorphic functions on $\R$, and $\R^+$,
respectively.; \item $C_0(J,\X):= \{ f\in BC(J,\X): \ \lim_{t\to
\infty} f(t)=0\}$; \item $AAP(\R^+,\X ) := C_0(\R^+,\X) \oplus
AP(\X)$; \item If $A$ is a linear operator on $\X$, then the
operator of multiplication by $A$ on $BC(J,\X)$, denoted by ${\cal
A}$, is defined on $D({\cal A}):= \{ g\in BC(J,\X):\ g(t)\in D(A),
\ \mbox{for all} \ t\in J, Ag(\cdot )\in BC(J,\X)\}$, by ${\cal
A}g=Ag(\cdot )$ for each $g\in D({\cal A})$.
\end{enumerate}
 In this paper, by almost periodic functions we mean the ones
in the sense of Bohr (for the precise definition and properties
see e.g. \cite{levzhi}), and by almost automorphic functions we
mean the ones in the sense of Bochner (for the precise definition
and properties see e.g. \cite{boc,ngu,sheyi,vee,yi}).

\section{A Spectral Theory of Bounded Functions on The Line}
\subsection{Reduced spectrum}
Let us introduce some operators and discuss the relations between
their resolvent sets and spectra of a bounded function.
\begin{definition}
$D({\cal D})$ is defined to be the set of all differentiable
functions $f\in BC(\R,\X)$ such that $f'\in BC(\R,\X)$. The
operator ${\cal D}$ is defined by ${\cal D}f=f'$ whenever $f\in
D({\cal D})$.
\end{definition}
Some elementary properties of the operator $\cD$ are summarized in
the following lemma.
\begin{lemma}\label{lem 1}
${\cal D}$ is a closed operator on $BC(\R,\X)$ with $\sigma ({\cal
D}) = i\R$. Moreover, for each $\xi\in\R$, $Re\lambda\not= 0$,
\begin{eqnarray}\label{re-1}
R(\lambda,{\cal D})f (\xi )&=& \begin{cases}
\begin{array}{ll}
\int^  {\infty}_0  e^{-\lambda \eta }f(\xi +\eta )d\eta
&(\mbox{if}\ Re\lambda > 0)\\ \\
-\int _ {-\infty}^0  e^{-\lambda \eta}f(\xi +\eta )d\eta     &
(\mbox{if } \ Re\lambda < 0).
\end{array} \end{cases}
\end{eqnarray}
\end{lemma}
\begin{proof}
Since for each $\lambda\in\C$ such that $Re\lambda \not= 0$, the
differential equation
\begin{equation}
x'(t)-\lambda x(t) =0, \quad x(t)\in \C
\end{equation}
has an exponential dichotomy, the non-homogeneous equation
\begin{equation}
x'(t)-\lambda x(t) =f(t), \quad x(t)\in \C
\end{equation}
has a unique solution $x_{f,\lambda}\in BC(\R,\X)$ for each given
$f\in BC(\R,\X)$. It is well known in the theory of ODE that
$x_{f,\lambda}$ is determined by the Green operator, that is,
\begin{eqnarray}\label{re-0}
-x_{f,\lambda }(\xi )  &=&
\begin{cases}
\begin{array}{ll}
\int^  {\infty}_\xi  e^{\lambda (\xi -t)}f(t)dt     &(\mbox{if}\
Re\lambda > 0)\\ \\
-\int _ {-\infty}^\xi  e^{\lambda (\xi -t)}f(t)dt    & (\mbox{if }
\ Re\lambda < 0)
 .
\end{array} \end{cases} \nonumber\\
&=&  \begin{cases}
\begin{array}{ll}
\int^  {\infty}_0  e^{-\lambda \eta }f(\xi +\eta )d\eta
&(\mbox{if}\ Re\lambda > 0)\\ \\
-\int _ {-\infty}^0  e^{-\lambda \eta}f(\xi +\eta )d\eta     &
(\mbox{if } \ Re\lambda < 0).
\end{array} \end{cases}
\end{eqnarray}
Therefore, $-x_{f,\lambda}$ is a bounded linear operator in $f$
acting in $BC(\R,\X)$, that is, $R(\lambda,{\cal D})$ exists and
 $R(\lambda,{\cal D})f=-x_{f,\lambda}$.

\medskip
Next, note that for each $i\xi$, where $\xi\in\R$, the function
$f_\xi (t):= e^{i\xi t}x$, where $x\in \X$ is a non-zero element,
is an eigenvector of $\cD$. Thus, $\sigma (\cD)=i\R$. As a
consequence, $\rho (\cD)\not= \emptyset$, so $\cD$ is closed.
Formula (\ref{re-1}) follows from (\ref{re-0}).
\end{proof}

Let $\F$ be a closed subspace of $BC(\R,\X)$ that satisfies the
following condition:
\begin{definition}
 A closed subspace $\F$ of $BC(\R,\X)$ is said to
satisfy {\it Condition F} if either it is trivial, or it satisfies
the following condition:
\begin{enumerate}
\item It contains all constant functions; \item If $f\in \F$,
then, for each $\xi\in R$, the function $f_\xi$ defined as $f_\xi
(t)=e^{i\xi t}f(t)$, for all $t\in \R$, belongs to $\F$; \item For
each $\lambda\in \C$ such that $Re\lambda\not= 0$, one has
\begin{equation}
R(\lambda ,\cD)  \F \subset \F ;
\end{equation}
\item If $f\in\F$ and $F$ is a bounded primitive of $f$, then
$F\in\F$; \item For each $B\in L(\X)$ and $f\in {\cal F}$, the
function $Bf(\cdot )$ is in $\cal F$.
\end{enumerate}
\end{definition}
\begin{remark}
\begin{enumerate}
\item In the case when $\F$ is a subspace of $BUC(\R,\X)$, condition (iii) follows from the translation invariance of $\F$. In fact, this follows from the representation
$$
R(\lambda ,\cD) g =\int^\infty _0 e^{-\lambda t}S(t)gdt, \quad Re\lambda >0 , t\in\R ,
$$
or
$$
R(\lambda ,\cD) g =-\int^\infty _0 e^{\lambda t}S(-t)gdt, \quad Re\lambda <0 , t\in\R .
$$
\item Condition iv) in the case $\F=AP(\X)$
is the validity of the well known Bohl-Bohr-Kadets Theorem (see
e.g. \cite{arebathieneu,levzhi}). For example, condition (iv) holds in this case if $\X$ does not contain any subspace isomorphic to $c_0$. For related concepts see \cite[Definition 3.1]{basgun}\footnote{Compared with the conditions listed in \cite[(3.1)]{basgun2} to define the concept of reduced spectrum our conditions are different. The uniform-closedness of $\cal A$ as a condition of \cite[Definition 3.1]{basgun} seems to be very strict. In fact, it can be easily checked that such a space ${\cal A}$ cannot be $BUC(\R,\X), BC(\R,\X)$}
\end{enumerate}
\end{remark}

\begin{example}\rm
Another example of such a function space $\F$ is the space of all
almost automorphic functions $AA(\X)$. It is known (see e.g.
\cite{minnaingu}) that $AA(\X) \not\subset BUC(\R,\X)$. Therefore,
the condition iii) does not follow from the translation invariance
of $AA(\X)$. This condition can be checked directly. A particular
case of the main result in \cite[Theorem 3.2]{minnaingu} yields
that the condition iv) is fulfilled for this function space if
$\X$ do not contain any subspace isomorphic to $c_0$.
\end{example}

Consider the quotient space $\Y:= BC(\R,\X)/ \F$, where $\F$ is a
given closed subspace of $BC(\R,\X)$ that satisfies Condition F.
Every element of this quotient space is a class of functions in
$BC(\R,\X)$ that is denoted by $\tilde{f}$, where $f\in
BC(\R,\X)$.
\begin{definition}
$D(\tilde{{\cal D}})$ is defined to be the subset of $\Y$
consisting of all classes of functions that contain elements of
$D(\cD)$. The operator $\tD$ is defined by $\tD \tf=\tilde{f'} $
whenever $\tf\in D(\tD)$ containing $f\in D(\cD)$.
\end{definition}
Some elementary properties of $\tD$ are summarized in the
following
\begin{lemma}\label{lem 2}
$\tD$ is a closed linear operator on $\Y$ with $\sigma (\tD)
\subset \sigma (\cD)=i\R $. Moreover, for $Re\lambda\not= 0$,
\begin{equation}\label{2}
\| R(\lambda ,\tD)\| \le \| R(\lambda ,\cD)\| \le
\frac{1}{|Re\lambda |} .
\end{equation}
\end{lemma}
\begin{proof} Take any $\lambda\in\rho (\cD)$, that is, $Re\lambda \not= 0$.
We will show that $\lambda \in \rho (\tD)$. In fact, for any
$f\in BC(\R,\X)$, since $R(\lambda ,\cD)\F \subset \F$, we have
$R(\lambda ,\cD)\tf$ is contained in $\tg$ defined as the class
containing $R(\lambda ,\cD)f$. So, the equation
\begin{eqnarray}\label{1}
\lambda \tg - \tD \tg =\tf
\end{eqnarray}
has at least one solution as the class containing $R(\lambda
,\cD)f$ for each given $\tf\in\Y$. Moreover, we have
\begin{eqnarray*}
\| \tg \|_\Y &:=& \inf_{h\in \F} \| R(\lambda , \cD ) f + h\|
\\
&\le& \inf_{k\in \F} \| R(\lambda , \cD ) (f + k)\| \\
&\le& \| R(\lambda , \cD )\| \inf_{k\in \F} \|  f + k\| \\
&=& \| R(\lambda , \cD )\| \cdot \|  f \|_\Y .
\end{eqnarray*}
Now we show that (\ref{1}) has no more than one solution. Indeed,
it is equivalent to show that the homogeneous equation has no
solution other than the zero solution. In fact, if $\lambda \tg
-\tD \tg =0$, then, assuming $\tg$ contains $g\in D(\cD)$, we have
$$
\lambda g(t)-g'(t) =h(t), \quad \ \mbox{for all} \ t\in \R,
$$
where $h$ is a function in $\F$. However, since $h\in\F$ from the
condition iii) of Condition F, $g$ must be in $\F$, so $\tg =0$.

\medskip
Summing up all we have done above shows that $\rho(\tD) \supset
\rho (\cD) = \C \backslash i\R  \not= \emptyset .$ As a
consequence, $\tD$ is closed. The estimate (\ref{2}) follows from
the above estimate of $\|R(\lambda ,\tD)\|$ in terms of $\|
R(\lambda,\cD)\|$, and in turn, an estimate of $\| R(\lambda
,\cD)\|$ from (\ref{re-1})
\end{proof}

We are ready to define the concept of reduced spectrum of a
bounded function.
\begin{definition}\label{def 1}
Let $\F$ be a closed subspace function of $BC(\R,\X)$ that
satisfies Condition F, and let $f\in BC(\R,\X)$. Then, the {\it
reduced spectrum of $f$ with respect to $\F$}, denoted by
$sp_{\F}(f)$,  is defined to be the set of all reals $\xi$ such
that the complex function $R(\lambda ,\tD)f$, as a function of
$\lambda \in \C \backslash i\R$, has no analytic extension to any
neighborhood of $i\xi$ in $\C$. If $\F$ is trivial, we use the
notation $sp(f)$ instead of $sp_0(f)$, and call it simply {\it
spectrum} of $f$.
\end{definition}
\begin{remark}
If $f\in BUC(\R,\X)$, and $\F$ is a subspace of $BUC(\R,\X)$ that
satisfies Condition F, then $sp_{\F}(f)$ is the reduced spectrum
of $f$ as defined in \cite{arebat,arebathieneu,bas,bask,ruevu}.
\end{remark}
\begin{proposition}\label{pro 1}
Let $\F$ be a closed subspace of $BC(\R,\X)$ that satisfies
Condition F. Then the following assertions hold:
\begin{enumerate}
\item $sp_\F (f)$ is closed, and $sp_\F(f)=sp_\F(g)$ for all $g\in
\tf$; \item $sp_\F (f) = sp_\F (f(\cdot +c ))$, for any $c \in\R$,
$f\in BC(\R,\X)$; \item $sp_\F(Af(\cdot ))\subset sp_\F (f)$ for
each $f\in BC(\R,\X)$, $A\in L(\X)$; \item $sp_\F (f +g)\subset
sp_\F(f)\cup sp_\F(g)$ for all $f,g\in BC(\R,\X)$; \item $sp_\F
(f) \subset \Lambda $ if there are $f_n\in BC(\R,\X)$, $n\in\N$,
such that $\tf_n \to \tf\in \Y:= BC(\R,\X)/\F$, $sp_\F (f_n)
\subset \Lambda $ for all $n$, where $\Lambda$ is a closed subset
of $\R$; \item $sp_\F (f) \subset sp(f)$.
\end{enumerate}
\end{proposition}
\begin{proof}
Properties i-iv) follows immediately from the definition. We now
prove v). Let $\rho_0 \not\in \Lambda$. Since $\Lambda$ is closed,
there is a positive constant $r< dist (\rho_0 ,\Lambda )$. As in
the proof of \cite[Theorem 0.8, p. 21]{pru} or by \cite[Lemma
4.6.6, p. 295]{arebathieneu} we can prove that since $R(\lambda
,\tD)\tf_n$ is extendable to all $\bar B_r(i\rho_0)$, and
\begin{equation}
\| R(\lambda ,\tD)\tf_n\| \le \frac{2\| \tf\|}{| Re \lambda |},
\quad \ \mbox{for all} \ \lambda \in \bar B_r(i\rho_0)
\end{equation}
for sufficiently large $n\ge N$,  one has
\begin{equation}
\|R(\lambda ,\tD)\tf_n\| \le \frac{4\| \tf\|}{3r}, \quad \
\mbox{for all} \ \lambda \in \bar B_r(i\rho_0), n\ge N.
\end{equation}
Obviously, for every fixed $\lambda$ such that $Re \lambda \not=
0$ we have $R(\lambda ,\tD)\tf_n \to  R(\lambda ,\tD)\tf$. Now
applying Vitali's theorem to the sequence of complex functions $\{
R(\lambda ,\tD)\tf_n\}$ we see that $R(\lambda ,\tD)\tf_n$ is
convergent uniformly on $B_r(i\rho_0)$ to $R(\lambda ,\tD)\tf$.
This yields that $R(\lambda ,\tD)\tf$ is holomorphic on
$B_r(i\rho_0)$, that is, $\rho_0 \not\in sp_\F (f)$.

\medskip For vi) it is obvious since the canonical projection on
the quotient space $BC(\R,\X) \to BC(\R,\X)/\F$ is continuous and
$p(R(\lambda ,\cD )f) = R(\lambda ,\tD )\tf$ for each $f\in
BC(\R,\X)$ and $Re\lambda \not= 0$.
\end{proof}
\begin{corollary}\label{cor 2}
Let $\F$ be a closed subspace of $BC(\R,\X)$ that satisfies
Condition F, and let $\Lambda$ be a closed subset of $\R$. Then
the function space
\begin{equation}
\Lambda_\F (\X) := \{ \tf\in  \Y := BC(\R,\X)/\F: \ sp_\F (\tf)
\subset \Lambda \}
\end{equation}
is a closed subspace of $\Y$ that satisfies
$$
R(\eta ,\tD)\Lambda_\F (\X)\subset \Lambda_\F (\X),\quad \
\mbox{for all} \ Re\eta \not= 0 .
$$
\end{corollary}
\begin{proof}
The first assertion of the corollary follows from Properties i-v)
of Proposition \ref{pro 1}. The last one follows from the note
that for $Re\lambda\not= 0$ and $Re\eta \not=0$,
$$
R(\lambda ,\tD )R(\eta ,\tD )\tf = R(\eta ,\tD)R(\lambda,\tD) \tf
.
$$
\end{proof}

Let $\F$ be a closed subspace of $BC(\R,\X)$ that satisfies
Condition F, and let $\Lambda$ be a closed subset of $\R$. Then,
we define an operator $\tD_\Lambda$ on $\Lambda_\F(\X)$ that is
the part of $\tD$ on $\Lambda_\F(\X)$, that is,
$$
D(\tD_\Lambda ):= \{ \tf \in D(\tD)\cap \Lambda_\F(\X):\ \tD \tf
\in \Lambda_\F (\X) \}
$$
and $\tD_\Lambda \tf = \tD\tf$, whenever $\tf\in D(\tD_\Lambda )$.
If $\F$ is trivial, we will use the notation $\cD_\Lambda$ instead
of $\tD_\Lambda$.

\begin{theorem} \label{the spe of tD Lambda}
Let $\F$ be a closed nontrivial subspace of $BC(\R,\X)$ that
satisfies Condition F. Then,
\begin{equation}\label{spe of tD Lambda}
\sigma (\tD_\Lambda )\subset i\Lambda .
\end{equation}
\end{theorem}
\begin{proof}
This is equivalent to show that every $\beta \in \R\backslash
\Lambda$ is in $\rho ({\cal D}_\Lambda )$, that is, the following
equation
\begin{equation}\label{3}
i\beta \tg-\tg'= \tf
\end{equation}
is solvable uniquely in $\Lambda_\F (\X)$ for every given
$\tf\in\Lambda_\F (\X)$. First, we show that (\ref{3}) has at most
one solution, or equivalently, the homogeneous equation $i\beta
\tg-\tg'=0$ has zero as the unique solution in $\Lambda_\F (\X)$.
In fact, if $\tg$ is a solution of this homogeneous equation, then
$$
i\beta g- g'= h \in \F .
$$
By the Variation-of-Constants Formula this equation is equivalent
to the following
$$
g(t)=e^{i\beta t}g(0) -\int^t_0 e^{i\beta (t-\xi )}h(\xi )d\xi
,\quad \ \mbox{for all} \ t\in\R .
$$
Therefore,
$$
e^{-i\beta t}g(t)=g(0) -\int^t_0 e^{-i\beta \xi }h(\xi )d\xi
,\quad \ \mbox{for all} \ t\in\R .
$$
Since $h\in\F$ and $F$ satisfies Condition F, we see that the
function $\R\ni \xi \mapsto e^{-i\beta \xi }h(\xi )$ belongs to
$\F$. Moreover, the function $\R\ni t\mapsto g(0)-e^{-i\beta
t}g(t)$ is a bounded primitive of the previous one, so by the item
i) of Condition F, the function $\R\ni t\mapsto e^{-i\beta t}g(t)$
belongs to $\F$. Therefore, $g$ also belongs to $\F$, that is,
$\tg=0$.

\medskip Next, we show that (\ref{3}) has at least one solution.
For every $Re\lambda \not= 0$ the equation
$$
\lambda y-y'=\tf
$$
has a unique solution $\tg_\lambda  = R(\lambda,{\tD})\tf$, that
is in $\Lambda_\F(\X)$ by Corollary \ref{cor 2}. Since
$R(\eta,{\tD}) \tf$ is analytic around $i\beta$, $\lim_{\lambda\to
i\beta}R(\lambda ,\tD)\tf $ exists as an element, say, $\tg\in
\Lambda_\F (\X)$. Now we show that $\tg$ is a solution of
(\ref{3}). Indeed, since
\begin{eqnarray*}
(i\beta -{\tD})R(\lambda ,\tD)\tf &=& ((i\beta -\lambda) + (\lambda -{\tD}))R(\lambda ,\tD)\tf \\
&=& (i\beta -\lambda )R(\lambda ,\tD)\tf +(\lambda -{\tD})R(\lambda ,\tD)\tf\\
&=& (i\beta -\lambda )R(\lambda ,\tD)\tf +\tf,
\end{eqnarray*}
and $R(\lambda ,\tD)\tf$ has an analytic extension around
$i\beta$, we have
\begin{equation}
\lim_{\lambda \to i\beta } (i\beta -{\tD})R(\lambda ,\tD)\tf =f.
\end{equation}
By the closedness of the operator $(i\beta -{\tD})$, we come up
with $\tg $ being in the domain of $i\beta -{\tD}$ and $(i\beta
-{\tD})\tg =\tf$.
\end{proof}
\begin{remark}
When $\F = \{ 0\}$, it is proved in \cite{liunguminvu} that
$\sigma (\cD _\Lambda ) =i\Lambda $. We refer the reader to
\cite{dav} for more related results concerned with the case $f\in
BUC(\R,\X)$. Results of this type can be used to study the
existence and uniqueness of bounded solutions to non-homogeneous
equations (see
\cite{murnaimin,murnaimin2,diangumin,liunguminvu,naiminshi}).
\end{remark}

\subsection{Coincidence of the notions of spectrum}
We first recall some concepts.
\begin{definition}\label{def 2}
The Carleman spectrum of $f\in BC(\R,\X)$ is defined to be the set
of all reals $\xi$ such that the Carleman transform
\begin{equation}
\hat{f}(\lambda) := \begin{cases}
\begin{array}{ll}
\int^  {\infty}_0  e^{-\lambda \eta }f( \eta )d\eta
&(\mbox{if}\ Re\lambda > 0)\\ \\
-\int _ {-\infty}^0  e^{-\lambda \eta}f( \eta )d\eta     &
(\mbox{if } \ Re\lambda < 0),
\end{array} \end{cases}
\end{equation}
as a complex function of $\lambda$, has no analytic extension to
any neighborhood of $i\xi$ .
\end{definition}
From the definition of Carleman spectrum of $f\in BC(\R,\X)$, that
will be denoted by $sp_c(f)$, it is clear that $sp_c(f) \subset
sp(f)$.

\begin{definition}\label{def 3}
Let $f\in BC(\R,\X)$. The Beurling spectrum of $f$, that is
denoted by $sp_b(f)$, is defined to be the following set
\begin{equation}
sp_b(f):= \{ \xi \in \R:\ \ \mbox{for all} \ \epsilon >0, \exists
\phi\in L^1(\R), \ supp (\tilde{\phi} ) \subset (\xi-\epsilon, \xi
+\epsilon ), \phi* f \not= 0 \} ,
\end{equation}
where
$$
\tilde{\phi } (\eta ) := \int^\infty_{-\infty} e^{i\eta t}\phi
(t)dt ,\quad t\in\R ,
$$
is the Fourier transform of $\phi$, and
$$
\phi*f (t) := \int^\infty_{-\infty} \phi (s) f(t-s )ds ,\quad \
\mbox{for all} \ t\in \R .
$$
\end{definition}
\begin{proposition}
Let $f\in BC(\R,\X)$. Then
\begin{equation}
sp(f)=sp_c(f)=sp_b(f) .
\end{equation}
\end{proposition}
\begin{proof}
For the proof of $sp_c(f)=sp_b(f)$ see \cite[Proposition 4.8.4, p.
321]{arebathieneu} and \cite{pru}. For the identity
$sp(f)=sp_c(f)=sp_b(f)$ see \cite{liunguminvu}.
\end{proof}

\subsection{Loomis Theorem}
The main result we will prove in the section is of the Loomis
Theorem type for general classes of functions. Before doing so, we
need some preparatory results.
The following lemma is known.
\begin{lemma}\label{lem prep gel}
Let $f(z)$ be a complex function taking values in a Banach space $\X$ and be holomorphic in $\C \backslash i\R$ such that there is a positive number $M$ independent of $z$ for which
\begin{equation}
\| f(z)\| \le \frac{M}{|Re\ z|}, \quad \ \mbox{for
all} \ \ Re z \not= 0 .
\end{equation}
Assume further that $i\xi\in i\R$ is an isolated singular point of $f(z)$ at which the Laurent expansion is of the form
\begin{equation}\label{laurent exp}
f(z) = \sum_{n=-\infty}^{\infty} a_n (z -i\xi )^n,
\end{equation}
where
\begin{equation}\label{laurent coeff}
a_n =\frac{1}{2\pi i} \int_{| z-i\xi |=r} \frac{f(z)dz }{(z-i\xi
)^{n+1}},\quad n\in \Z .
\end{equation}
Then,
\begin{equation}\label{an}
\| r^2a_{-(n+1)} + a_{-(n+3)}\| \le 2Mr^{n+2}, \quad n\in \Z .
\end{equation}
\end{lemma}
\begin{proof}
Note that the below proof can be found in \cite[Lemma 4.6.6]{arebathieneu}, \cite[Chap. 0]{pru}. For the reader's convenience we reproduce it here.

\medskip
For each $n\in \Z$ and
$0<r<\delta_0$, where $\delta_0$ is some positive number, we have
\begin{eqnarray*}
&&\| \frac{1}{2\pi i} \int_{|z-i\xi |=r} (z-i\xi )^n \left(
1+\frac{(z-i\xi )^2}{r^2}\right) f(z)dz \| \\
&& \hspace{2cm} \le \frac{1}{2\pi} \int_{|z-i\xi|=r} |(z-i\xi )^n
\left( 1+\frac{(z-i\xi )^2}{r^2}\right)|\cdot \| f(z)\|\cdot |dz | .
\end{eqnarray*}
A simple computation shows that since $|z-i\xi |=r$, one has
\begin{equation}\label{n15}
|(z-i\xi )^n \left( 1+\frac{(z-i\xi )^2}{r^2}\right)| = 2r^{n-1} | Re\
z | .
\end{equation}
Therefore,
\begin{eqnarray}
\| \frac{1}{2\pi i} \int_{|z-i\xi |=r} (z-i\xi )^n \left(
1+\frac{(z-i\xi )^2}{r^2}\right) f(z)dz \| &\le& \frac{1}{2\pi}
\int_{|z-i\xi|=r}2r^{n-1} | Re\ z | \frac{M}{|Re \   z|} \cdot |dz |
\nonumber \\
&=& \frac{2Mr^{n-1}}{2\pi} \int_{|z-i\xi|=r}  |dz | \nonumber \\
&=&{2Mr^{n}} .\label{15}
\end{eqnarray}
Consider the Laurent expansion (\ref{laurent exp}).
From (\ref{15}) it follows that for all $n\in\Z$,
\begin{eqnarray}
\| a_{-(n+1)} + r^{-2} a_{-(n+3)}\| &=& \| \frac{1}{2\pi
i}\int_{|z-i\xi |=r} (z-i\xi )^n f(z)dz
 + \frac{1}{2\pi i} \int_{|z-i\xi |=r} \frac{(z-i\xi )^{n+2}}{r^2}
f(z)dz \| \nonumber\\
&=&  \| \frac{1}{2\pi i} \int_{|z-i\xi |=r} (z-i\xi )^n \left(
1+\frac{(z-i\xi )^2}{r^2}\right) f(z)dz \|  \nonumber\\
&\le& 2Mr^{n}.\label{n21}
\end{eqnarray}
Multiplying both sides by $r^2$ gives (\ref{an}). The lemma is proven.
\end{proof}

Below we offer a simple proof of an extension of the Gelfand Theorem for groups of isometries (see e.g. \cite{arebat}, \cite{arebathieneu}, \cite{nee} for more information about this theorem and applications). For the reduced spectrum of not necessarily uniformly continuous functions that is presented below this extension will play exactly the role of the Gelfand Theorem for groups of isometries in the study of the reduced spectrum of uniformly continuous functions in \cite{arebat,arebathieneu}.

\begin{theorem}\label{the gel}
Let $A$ be a closed linear operator on a Banach space $\X$ such
that
\begin{enumerate}
\item $\sigma (A) \subset i\R $; \item For some
$\lambda$-independent positive number $M$, the following condition
holds
\begin{equation}
\| R(\lambda ,A)\| \le \frac{M}{|Re\ \lambda|}, \quad \ \mbox{for
all} \ \ Re \lambda \not= 0 .
\end{equation}
\end{enumerate}
Then, the following assertions hold
\begin{enumerate}
\item If $\lambda_0=i\xi \in i\R$ is an isolated point of $\sigma (A)$,
then it is an eigenvalue of $\sigma (A)$; \item If $\X$ is non-trivial, then $\sigma (A) \not=
\emptyset$;
\item If $\sigma (A) =\{ 0\} $, then $A=0$.
\end{enumerate}
\end{theorem}
\begin{proof}
i) Set $f(\lambda )= R(\lambda ,A)$. Since $i\xi \in i\R$ is
an isolated point in $\sigma (A)$, it is an isolated singular point of $f(\lambda )$, so by Lemma \ref{lem prep gel}
\begin{equation}\label{an 2}
\| r^2a_{-(n+1)} + a_{-(n+3)}\| \le 2Mr^{n+2}, \quad n\in \Z .
\end{equation}
 Letting
$r$ tend to $0$ in (\ref{an}), we come up with $a_{-k}=0$ for all
$k \ge 2$. This shows that $i\xi$ is a pole of first order of the
resolvent $f(\lambda ) :=R(\lambda ,A)$. And hence, by a well
known result in Functional Analysis (see e.g. \cite[Theorem 5.8 A,
p. 306]{tay}, or, \cite[Theorem 3, p. 229]{yos}), $i\xi$ is an
eigenvalue of the operator $A$. So, the first assertion is proved.

\medskip
ii) Next, suppose that $\rho(A)=\C$. Consider the Laurent
expansion (\ref{laurent exp}) of $f(\lambda )=R(\lambda,A)$ at
$\lambda =i\xi$. Since $R(\lambda ,A)$ is analytic everywhere, $a_n=0$ for all $n\le -1$. Note that
the formula (\ref{n21}) is still valid for this case, and can be
re-written in the form
\begin{equation}\label{an positive}
\| a_{k-1} + r^{-2}a_{k-3}\| \le 2M\frac{1}{r^{k}}, \quad k\in \Z .
\end{equation}
Letting $r$ tend to infinity we have $ a_n=0$ whenever $n\ge 0$, so $R(\lambda ,A)=0$.
This is impossible if $\X$ is non-trivial. This contradiction
proves ii).

\medskip
iii) By (\ref{an}) and (\ref{an positive}) it is easy to see that all
$a_n=0$ with $n\not = -1$. So,
$$
R(\lambda ,A)= \frac{a_{-1}}{\lambda } , \quad \ \mbox{for
all} \ \ \lambda \not = 0 .
$$
We have
$$
I=\left(\lambda - A\right)R(\lambda ,A)=\left(\lambda -
A\right)\left(\frac{a_{-1}}{\lambda }  \right)= a_{-1}
  -\frac{Aa_{-1}}{\lambda}  , \quad \ \mbox{for
all} \ \ \lambda\not= 0.
$$
Letting $\lambda$ tend to infinity we can show
that $a_{-1}=I$, and thus, $R(\lambda
,A)=I/\lambda$ for $\lambda\not= 0$. However, this yields $I=
I-A/\lambda $ for all $\lambda \not= 0,$ so, $A=0$.
\end{proof}
\begin{remark}
When $A$ is the generator of a bounded $C_0$-group
$(T(t))_{t\in\R}$, it satisfies the assumptions of these lemmas.
In fact, in this case, since
$$
R(\lambda ,A)x=\int^\infty _0 e^{-\lambda t}T(t) x dt, \quad
x\in\X, \ Re\lambda >0,
$$
and
$$
R(\lambda ,A)x= -\int^{\infty} _0 e^{\lambda t}T(-t) x dt, \quad
x\in\X, \ Re\lambda <0,
$$
 we have
$$
\| R(\lambda ,A) \| \le \frac{M}{|Re\lambda | }   ,
$$
where $M:= \sup_{t\in\R} \| T(t)\|$. Therefore, the above lemmas
extend the well-known Gelfand's theorem for $C_0$-groups of
isometries (see e.g. \cite[Corollaries 4.4.8 \&
4.4.9]{arebathieneu}). This will be a key point in the framework
for the next results.
\end{remark}

\begin{corollary}
Let $\F$ be a subspace of $BC(\R,\X)$ that satisfies Condition F,
and let $f\in BC(\R,\X)$. Then $sp_\F(f) =\emptyset$ if and only
if $f\in \F$.
\end{corollary}
\begin{proof}
If $f\in \F$, then the assertion is obvious. Conversely, let
$sp_\F(t)=\emptyset$. Set $\Lambda =\emptyset$. If $\Lambda_\F
(\X)$ is non-trivial, then, by Theorem \ref{the spe of tD Lambda},
$\sigma (\tD) \subset i\Lambda$, so $\sigma (\tD)=\emptyset$. This
contradicts Theorem \ref{the gel}. Therefore, $\Lambda_\F (\X)$ is
trivial, and $f\in \F$.
\end{proof}

The following corollary is well known in the spectral theory of
functions (see e.g. \cite{pru,arebathieneu}). However, we will
restate it and give a proof based on the our approach to the
spectrum.
\begin{corollary}\label{cor 2.19}
Let $f\in BC(\R,\X)$. Then $sp (f) =\{ \xi_1, \xi_2, \cdots
,\xi_N\}$ if and only if $f$ is of the form $\R \ni t \mapsto
\sum_{k=1}^Na_ke^{i\xi_kt}\in\X$, where $0\not= a_k\in \X $.
\end{corollary}
\begin{proof}
The necessity is obvious. Now we show the sufficiency. Set
$\Lambda =sp (f)$. By the Riesz decomposition of closed operators
(see e.g. \cite[Chap. 4]{engnag}) we can decompose $\Lambda (\X)$
as $\Lambda  (\X)= \Lambda^1 \oplus \cdots \oplus \Lambda ^{N
}\oplus \Lambda^{N+1}$, where the spectrum of the restriction of
$\cD_\Lambda$ to $\Lambda^k$, a closed subspace of $\Lambda_\F
(\X)$,
 is contained in $\{ \xi_k\}$ for all $k=1,\cdots ,k$, and the spectrum of
 the restriction of $\cD_\Lambda$ to $\Lambda^{N+1}$ is
 empty. Moreover, the restrictions of $\cD_\Lambda$ to
 $\Lambda^k$, $k=1,\cdots, N$ are bounded.
By Theorem \ref{the gel} $\Lambda^{N+1}$ must be trivial. Therefore,
it suffices to show that if $sp  (f) =\{ \xi \}$, then $f$ is
 of the form $ae^{i\xi t}$ with
$a\not= 0$. Without loss of generality we may assume that $\xi
=0$. By Theorem \ref{the spe of tD Lambda}, $ \sigma (\cD_\Lambda
) = \{ 0\} $, so, by Theorem \ref{the gel}, $\cD_\Lambda =0$. Since
$\cD_\Lambda$ is bounded and $D(\cD) =\Lambda (\X)$, so $f\in
D(\cD) =\Lambda (\X)$. And we have $f'=o$. This shows $f(t)=
const$ . The corollary is proved.
\end{proof}
\begin{corollary}
Let $f\in BC(\R,\X)$. If $\xi$ is an isolated point in $sp (f)$,
then $\xi \not\in sp_{AP(\X)}(f)$.
\end{corollary}
\begin{proof}
Set $\Lambda := sp(f)$. By Theorem \ref{the spe of tD Lambda},
$\Lambda (\X)$ can be decomposed as $\Lambda (\X)= \Lambda_1
\oplus \Lambda_2$, where the restriction of $\cD_\Lambda$ to
$\Lambda_1$ is bounded and has spectrum as $\{ \xi\}$, and the
restriction of $\cD_\Lambda$ to $\Lambda_2$ has the spectrum as
$sp(f)\backslash \{ \xi\}$. Therefore, $f= f_1+f_2$, where
$sp(f_1)=\{\xi\}$ and $sp(f_2) \subset sp(f)\backslash \{ \xi\}$.
By Corollary \ref{cor 2.19}, $f_1$ is of the form $f_1(t)
=ae^{i\xi t}$. Hence, $\tf = \tf_1=\tf_2 =\tf_2$. By Theorem
\ref{the spe of tD Lambda} $sp_{AP(\X)}(f_2) \subset sp(f_2)
\subset sp(f)\backslash \{ \xi\}$. Therefore, $\xi\not\in
sp_{AP(\X)}(f)$.
\end{proof}
The above corollary is known in \cite{aresch} with additional
assumption that $f\in BUC(\R,\X)$. An immediate consequence of
this lemma is the following:
\begin{corollary}
Let $f\in BC(\R,\X)$. If $sp (f)$ is discrete, then $f$ is almost
periodic.
\end{corollary}
\begin{remark}
This corollary with additional assumption on the uniform
continuity of $f$ has been known in \cite{aresch,arebathieneu},
and in a more abstract contexts in \cite{bask,beu,rei}.
\end{remark}

\bigskip
We are in a position to prove the following that is often referred
to as the Loomis Theorem, or, of Loomis Theorem type.
\begin{theorem}\label{the loomis}
Let $\F$ be a closed subspace of $BC(\R,\X)$ that satisfies
Condition F, and let $f\in BC(\R,\X)$ with countable $sp_\F(f)$.
Then, $f$ is in $\F$.
\end{theorem}
\begin{proof}
Let $\Lambda := sp_\F(f)$. We will show that $\Lambda_\F(\X)$ is
trivial. Suppose to the contrary that $\Lambda_\F(\X)$ is
non-trivial. Then, by Theorem \ref{spe of tD Lambda} and the
assumption, $\sigma (\tilde{\cal D}_\Lambda )$ is countable. By
Theorem \ref{the gel}, it is non-empty, so, since it has an isolated
point, it has an eigenvalue, say, $i\xi$, where $\xi\in\R$.
 This means that
there exists a non-zero $\tilde{g}\in D(\tilde{\cal D}_\Lambda)
\subset \Y$ such that
$$
\tilde{\cal D}\tilde{g}-i\xi \tilde{g}=0 .
$$
Therefore, the class $\tilde{g}$ contains a differentiable
function, say $g\in BC(\R,\X)$, and
$$
g'-i\xi g= h \in \F .
$$
Using the Variation-of-Constants Formula we have
$$
g(t) = e^{i\xi t}g(0) + \int^t_0 e^{i\xi (t-\eta)} h(\eta )d\eta ,
\quad \ \mbox{for all} \ t\in\R .
$$
Therefore,
$$
e^{-i\xi t}g(t) = g(0) + \int^t_0 e^{ -i\xi \eta } h(\eta )d\eta ,
\quad \ \mbox{for all} \ t\in\R .
$$
Since the function $\R\ni \eta \mapsto e^{-i\xi  \eta } h(\eta )$
is in $\F$, and its primitive $\R\ni t\mapsto e^{-i\xi t}g(t)$ is
bounded, by Condition F, the primitive $\R\ni t\mapsto e^{-i\xi
t}g(t)$ is in $\F$. Hence, $\tilde{g}=\tilde{0}$. This leads to a
contradiction proving that $\Lambda_\F(\X)$ is trivial.
\end{proof}
\begin{remark}
The above theorem has been stated and proved (see
\cite{arebat,arebathieneu,bas,basgun,baspry,bask,bask2}) under an
additional assumption on the uniform continuity of $f$. This
assumption is essential for the use of the techniques involving
the theory of $C_0$-groups.
\end{remark}

Some standard corollaries of Theorem \ref{the loomis} are as
follows:
\begin{corollary}\label{cor scalar}
Every scalar bounded and continuous function on $\R$ whose
spectrum is countable is almost periodic.
\end{corollary}
\begin{proof}
Let ${\cal F}:= AP(\R)$. By the Bohl-Bohr Theorem saying that
every bounded primitive of an almost periodic (scalar) function is
almost periodic, we can see that ${\cal F}$ satisfies Condition F.
Since $sp_{\cal F}(u)\subset sp(f)$, $sp_\F (f)$ is countable. So,
by Theorem \ref{the loomis}, $u\in {\cal F}$, that is, $u$ is
almost periodic.
\end{proof}
\begin{remark}
The above corollary seems to be new even in the scalar case. In
fact, in \cite{loo} Loomis proved the above corollary in a larger
context but with additional assumption on the uniform continuity
of $u$. The uniform continuity assumption is essential for the
techniques used in subsequent extensions (see
\cite{arebat,arebathieneu,bas,basgun,baspry,bask,bask2,ruevu}).
\end{remark}
\begin{corollary}\label{cor loomis 2}
Let $\X$ be a Banach space which does not contain any subspace
isomorphic to $c_0$. Then, every bounded and continuous function
with countable spectrum is almost periodic.
\end{corollary}
\begin{proof}
Let ${\cal F}:= AP(\X)$. By the Bohl-Bohr-Kadets Theorem, every
bounded primitive of an almost periodic function taking values in
a Banach space $\X$ not containing any subspace isomorphic to the
space of $c_0$, is almost periodic. So, the function space ${\cal
F}:= AP(\X)$ satisfies Condition F. Now by the same argument as in
the proof of the above corollary we can prove the corollary.
\end{proof}
\begin{remark}
The above corollary was first proved by Zhikov (see \cite{levzhi})
with additional assumption on the uniform continuity.
\end{remark}
Let us consider an example with ${\cal F}:= AA (\X)$, where $AA
(\X)$ denotes the space of all almost automorphic functions
introduced by Bochner. For the precise definition and properties
of these functions see e.g. \cite{boc,minnaingu,ngu,sheyi,yi,vee}.
As a special case, the main result in \cite{minnaingu} actually
says that if $\X$ does not contain any subspaces isomorphic to
$c_0$, then each bounded primitive of an $\X$-valued almost
automorphic function is almost automorphic. That is, ${\cal F}:=
AA (\X)$ satisfies Condition F in this case. Therefore, we arrive
at
\begin{corollary}\label{cor aa}
Let $\X$ be a Banach space which does not contain any subspaces
isomorphic to $c_0$, and let $f$ be in $BC(\R,\X)$ with countable
$sp_{AA(\X)}(f)$. Then, $f$ is in $AA(\X)$.
\end{corollary}
\begin{remark}
As $AA(\X) \not\subset BUC(\X)$, the above corollary seems to be
new.
\end{remark}

\bigskip
Before closing this subsection we would like to emphasize that we
can derive a version of Theorem \ref{the loomis} in which ${\cal
F}$ satisfies all conditions of Condition F except for the
condition iv) that is replaced by an ergodicity condition as
discussed in the next section. To avoid repeating the argument in
the next subsection we will state the ergodicity condition only
for the results for the functions on the half line. The reader can
easily adapt them to the entire real line case.

\section{Functions on the Half Line}
Let us consider differential equations of the form
\begin{equation}\label{hl}
\dot x(t) = \lambda x(t) +f(t), \quad t\in \R^+ ,
\end{equation}
where $f\in BC(\R^+,\X)$. If $Re\lambda >0$, the general solution
of (Eq. (\ref{hl}) is
\begin{equation}
x(t) = e^{\lambda t} x - \int^\infty _t e^{\lambda (t-s)} f(s)ds
,\quad x\in \X, \ t\in \R^+ .
\end{equation}
Therefore, the only bounded solution of (\ref{hl}) is
\begin{equation}\label{hl 1}
x_{\lambda ,f} (t):=- \int^\infty _t e^{\lambda (t-s)} f(s)ds , \quad
t\in \R^+ .
\end{equation}
On the other hand, if $Re\lambda <0$, the general solution of
(\ref{hl}) is
\begin{equation}\label{hl gs}
x(t) = e^{\lambda t} x + \int^t_0 e^{\lambda (t-s)} f(s)ds ,\quad
x\in \X, \ t\in \R^+ ,
\end{equation}
so, all solutions in this case are bounded, and all approach zero,
except for
\begin{equation}\label{hl b}
x_{\lambda ,f} =  \int^t_0 e^{\lambda (t-s)} f(s)ds ,\quad x\in
\X, \ t\in \R^+ .
\end{equation}

\medskip
Let us consider a function space ${\cal F} \subset BC(\R^+,\X)$
that satisfies the following Condition $F^+$:
\begin{definition}
 A function space ${\cal F} \subset BC(\R^+,\X)$ is
said to satisfy {\it Condition $F^+$} if
\begin{enumerate}
\item It is closed, and contains $C_0(\R^+,\X)$; \item If $g\in
{\cal F}$, then the function $\R^+\ni t\mapsto e^{i\xi t}g(t)\in
\X$ is in ${\cal F}$ for all $\xi\in\R$; \item For each $h\in
{\cal F}$, $Re\lambda >0, \ Re\eta <0$, the function $y(\cdot),
z(\cdot )$, defined as
\begin{equation}
y(t) = \int^\infty _t e^{\lambda (t-s)} h(s)ds, \ \ z(t) =
\int^t_0 e^{\eta (t-s)} h(s)ds ,\quad  t\in \R^+
\end{equation}
are in ${\cal F}$;  \item For each $B\in L(\X)$ and $f\in {\cal
F}$, the function $Bf(\cdot )$ is in $\cal F$.
\end{enumerate}
\end{definition}
As an example of a function space that satisfies Condition $F^+$,
we can take ${\cal F}=C_0(\R^+,\X)$. Another function space that
satisfies Condition $F^+$ is $AA(\R^+,\X)$, the space of all
restrictions to $\R^+$ of the $\X$-valued almost automorphic
functions. Note that $AA(\R^+,\X)$ contains non-uniformly
continuous functions, so it is not a subspace of $BUC(\R^+,\X)$
(see e.g. \cite{minnaingu}).

\bigskip
 In order to clarify the role condition (iii) of condition $F^+$ let us consider the differentiation operator ${\cD}$ on $BUC(\R^+,\X)$ to which we assume $\F$ belong. It is easy to see that for $Re\lambda >0$, $\lambda \in \rho ({\cD})$, and
$$
[R(\lambda ,{\cD}) f ](t)=-x_{\lambda ,f}(t)=\int_t^\infty e^{\lambda (t-s)}f(s)ds.
$$
Therefore, the first part of condition (iii) means that $R(\lambda ,{\cD})\F \subset \F$ (because $y=R(\lambda ,{\cD})f \in \F$ for each $f\in\F$). Since ${\cD}$ generates the translation semigroup, using the representation
$$
R(\lambda ,{\cD})f=\int^\infty _0 e^{-\lambda t}S(t)dt ,\quad Re\lambda>0,
$$
we can see that this first part of condition (iii) is satisfied if $\F$ is left invariant under the translation semigroup $S(t)$ in $BUC(\R,\X)$, that is, $S(t)\F \subset \F$ for all $t\ge 0$. The inverse is also true from the semigroup theory.

\medskip
Of course, in $BUC(\R^+,\X)$ the translation $S(t)$ is not invertible ($t>0$), and $\sigma (\cD )$ contains all complex numbers $z$ with $Rez<0$. However, the formulas (\ref{hl gs}) and (\ref{hl b}) give some insights into the structure of solutions of Eq. (\ref{hl}). "Within" an asymptotically stable solution $R(\lambda ,{\cD})f$ can be determined uniquely by (\ref{hl b}) even if $Re\lambda <0$. This leads to the idea of factoring all functions by asymptotically stable functions so that $Re\lambda <0$ belongs to the resolvent set of $\cD$ in the quotient space. And this will be all complex plane, but $i\R$. Our second part in condition (iii), that says that $z(t)$ is in $\F$, aims at realizing this idea.
This is crucial step for us to use the Gelfand Theorem (Theorem \ref{the gel}) to study the reduced spectrum of functions on the half line.

\medskip
In the approach to the reduced spectrum concept via the translation semigroup Arendt and Batty introduced the concept of biinvariance of $\F$ with respect to the translation semigroup $(S(t)_{t\ge 0}$ in $BUC(\R^+,\X)$, that is, the condition $S(t)\F =\F$ for all $t\ge 0$ (see e.g. \cite{arebat,arebathieneu}). This yields the surjectiveness of the isometries in the semigroup induced by this translation semigroup. And the Gelfand Theorem can be applied to study the reduced spectrum concept defined in this way.

\medskip
In the case $\F \subset BUC(\R^+,\X)$ since the induced differentiation operator $\cD$ in the quotient space mentioned above generates the induced translation semigroup we can show that our condition (iii) is equivalent to the biinvariance condition. The advantage of using condition (iii) is clear when it comes to $BC(\R^+,\X)$ in which the translation semigroup is not strongly continuous. 

\medskip
In \cite{basgun2} a concept reduced spectrum of not necessarily uniformly continuous functions on the whole line is also defined. Note that the uniform closedness condition (see \cite[Def. 3.1 and Theorem 4.3]{basgun2}) seems to be too restrictive.\footnote{As we have noted above in the theorem cited the authors obviously assume the non-emptyness of the spectrum.} Moreover, the results need to be adjusted to apply to equations on the half line.

\bigskip
Consider the quotient space $\Y := BC(\R^+,\X)/ {\cal F}$. We will
use $\tilde{\cD}$ to denote the operator induced by $\cD$ on $\Y$
which is defined as follows: The domain $D(\tD )$ is the set of
all classes that contains a differentiable function $g\in
BC(\R^+,\X)$ such that $g'\in BC(\R^+,\X)$; $\tD \tg  := \tg '$
for each $\tg \in D(\tD )$.

\bigskip
By (\ref{hl}) and (\ref{hl 1}), and the axiom iii) of Condition
$F^+$,
\begin{eqnarray}\label{resolvent of tilde D}
R(\lambda ,\tD )\tf  (t) &=& \begin{cases} \int^\infty _t
e^{\lambda
(t-s)} \tf (s)ds  ,\quad  Re\lambda > 0, \ t\in \R^+ ,\\ ~ \\
- \int^t_0 e^{\lambda (t-s)} \tf (s)ds  \quad Re\lambda <0, \ \
t\in \R^+.
\end{cases}\\
&=&\begin{cases} \int^\infty_0 e^{-\lambda \eta} \tf (t+\eta )d\eta ,\quad  Re\lambda > 0, \ t\in \R^+ ,\\ ~ \\
-\int^t_0 e^{\lambda \eta}\tf (t-\eta ) d\eta , \quad Re\lambda
<0, \ \ t\in \R^+.\nonumber
\end{cases}
\end{eqnarray}
\begin{lemma}\label{lem 4}
Under the above notations, the operator $\tilde{\cD}$ is a closed
operator with $\sigma (\tilde{\cD} ) \subset i\R$. Moreover, for
$Re\lambda \not= 0$,
\begin{equation}\label{hl2}
\| R(\lambda ,\tilde{\cD}) \| \le \frac{1}{| Re\lambda |} .
\end{equation}
\end{lemma}
\begin{proof}
By the above observations, the first assertion of the lemma is
obvious. Next, to show (\ref{hl2}) we can use axiom iii) of
Condition $F^+$. Therefore, by definition, for $Re\lambda
>0$, we have
\begin{eqnarray}
\| R(\lambda ,\tilde{\cD}) \tilde{f}\|_\Y &=& \inf_{g\in {\cal F}}
\sup_{t\in \R^+}\| g(t) + \int^\infty _t e^{\lambda
(t-s)} f(s) ds  \| \nonumber \\
&\le & \inf_{h\in {\cal F}} \sup_{t\in \R^+}\|  \int^\infty _t
e^{\lambda (t-s)} [f(s)+h(s)]  ds  \|\nonumber \\
&\le&\inf_{h\in {\cal F}}  \|  \int^\infty _t e^{Re \lambda
(t-s)}ds \| f+h\| \nonumber \\
&=& \|  \int^\infty _t e^{Re \lambda
(t-s)}ds \inf_{h\in {\cal F}} \| f+h\| \nonumber \\
&=& \frac{1}{|Re\lambda |} \| \tilde{f}\|_\Y .
\end{eqnarray}
Similarly, for $Re\lambda <0$ we can show that (\ref{hl2}) holds.
This proves the lemma.
\end{proof}
\begin{definition}
Let ${\cal F}$ be a function space that satisfies Condition $F^+$,
and let $f\in BC(\R^+,\X)$. Then the {\it reduced spectrum} of $f$
with respect to ${\cal F}$, denoted by $sp^+_{{\cal F}}(f)$, is
defined to be the set of all reals $\xi \in\R$ such that
$R(\lambda , \tilde{\cD})\tilde{f}$, as a complex function of
$\lambda$ in $\C \backslash i\R$, has no holomorphic extension to
any neighborhood of $i\xi$ in the complex plane.
\end{definition}
Since $sp^+_{{\cal F}}(f)$ is the same for all elements $f$ in a
class $\tg$, the use of the notation $sp^+_{{\cal F}}(\tg)$ makes
sense.
\begin{definition}
Let $\Lambda$ be a closed subset of the real line. Then,
\begin{equation}
\Lambda^e _{\cal F}(\X) := \{ \tf\in BC(\R^+,\X)/{\cal F}:\
sp_{{\cal F}}^+ (\tf) \subset \Lambda ,\ \lim_{\alpha \downarrow
0} \alpha R(\alpha +i\xi ,\tD )\tf =\tilde{0}\ \mbox{ for all $\xi
\in \Lambda $ }\} .
\end{equation}
\end{definition}
The property that $\lim_{\alpha \downarrow 0} \alpha R(\alpha
+i\xi , D ) f =g\in BC(\R,\X)$ is often referred to as the {\it
uniform ergodicity} of $f$ at $i\xi \in i\R$. For related concepts
of ergodicity and their equivalence to this one we refer the
reader to \cite{arebat2,arebathieneu,bas,delvu,ruevu}.

\medskip
Let us consider the restriction $\tD_\Lambda$ of $\tD$ to
$\Lambda^e _{\cal F}(\X)$.

\begin{theorem}\label{the spe dif ope hl}
Let $\Lambda$ be a closed subset of the real line. Then,
$\Lambda^e _{\cal F}(\X)$ is a closed subspace of $\Y:=
BC(\R^+,\X)/{\cal F}$, and
\begin{equation}\label{spe of tD hl}
\sigma (\tilde{\cD}_\Lambda ) \subset i\Lambda .
\end{equation}
\end{theorem}
\begin{proof}
To show the closedness of $\Lambda^e _{\cal F}(\X)$ we assume that
$\{\tf_n\}_{n=1}^\infty \in  \Lambda^e _{\cal F}(\X)$ such that
$\tf_n \to \tf\in \Y$ as $n\to \infty$. Using exactly the argument
in the proof of Proposition \ref{pro 1} we can easily show that
$sp_{{\cal F}}^+ (\tf)\subset \Lambda$. Next we will show that
\begin{equation}\label{erg}
\lim_{\alpha \downarrow 0} \alpha R(\alpha +i\xi ,\tD )\tf
=\tilde{0}.
\end{equation}
In fact, by the assumption, for each $\epsilon >0$ there is a
positive integer $N$ such that if $n\ge N$, then, $\| \tf_n
-\tf\|_\Y < \epsilon $. So, by (\ref{resolvent of tilde D}) for
each $\epsilon >0$ and sufficiently large $n$,
\begin{eqnarray*}
\limsup_{\alpha \downarrow 0} \alpha \| R(\alpha +i\xi ,\tD
)\tf\|_\Y &\le & \limsup_{\alpha \downarrow 0} \alpha \| R(\alpha
+i\xi ,\tD )\tf_n\|_\Y + \limsup_{\alpha \downarrow 0} \alpha \|
R(\alpha +i\xi ,\tD )(\tf_n -\tf )\|_\Y \\ &\le& 0+
\limsup_{\alpha \downarrow 0} \alpha \int^\infty e^{\alpha t}\|
\tf_n -\tf \|_\Y dt = \| \tf_n -\tf \|_\Y <  \epsilon .
\end{eqnarray*}
Since $\epsilon >0$ is arbitrary, this proves (\ref{erg}),
yielding the closedness of $\Lambda^e _{\cal F}(\X)$.

\medskip
Now we prove (\ref{spe of tD hl}), by solving $\tg\in \Lambda^e
_{\cal F}(\X)$ from the equation
\begin{equation}\label{100}
i\beta \tg - \tD_\Lambda \tg = \tf ,
\end{equation}
for each $\tf \in \Lambda^e _{\cal F}(\X)$, and $\beta \in\R$ such
that $\beta \not\in \Lambda $. First, we show the uniqueness.
Assume that $i\beta \tg - \tD_\Lambda \tg =0$.
 Since $(i\beta -\tilde{\cD})\tilde{g}=0$, we have
\begin{eqnarray}\label{3.13}
\tilde{g}&=& ((\alpha +i\beta ) -\tD )  R(\alpha +i\beta
,\tilde{\cD} )
\tilde{g}   \\
&=& \alpha   R(\alpha +i\beta ,\tilde{\cD} ) \tilde{g}
+ ( i\beta
 - \tilde{\cD}   )  R(\alpha +i\beta , \tilde{\cD} ) \tilde{g} \nonumber \\
&=&\alpha   R(\alpha +i\beta ,\tilde{\cD} ) \tilde{g} + R(\alpha
+i\beta , \tilde{\cD} )  ( i\beta
 - \tilde{\cD}   )   \tilde{g}  \nonumber \\
 &=&
 \alpha   R(\alpha +i\beta ,\tilde{\cD} ) \tilde{g}. \nonumber
\end{eqnarray}
On the other hand, since $\beta\not\in sp_{{\cal F}}^+ (\tf)$, the
function $R(\lambda ,\tilde{\cD} ) \tilde{g}$ has an analytic
extension to a neighborhood of $i\beta$. In particular,
$\lim_{\alpha\downarrow 0}R(\alpha +i\beta ,\tilde{\cD} )
\tilde{g}$ exists, so,
\begin{eqnarray*}
\tilde{g}&=&\lim_{\alpha\downarrow 0}\tilde{g}
 =  \lim_{\alpha\downarrow 0} \alpha   R(\alpha +i\beta ,\tilde{\cD} )
 \tilde{g}
 =  0.
\end{eqnarray*}
Now we prove the existence of a solution to (\ref{100}). Acting as
in the proof of Theorem \ref{the spe of tD Lambda} we can show
that $\tg := \lim_{\lambda\to i\beta}R(\lambda ,{\tD})\tf $ exists
as an element of $\Y$ such that $i\beta \tg -\tg ' =\tf$ and
$sp_{{\cal F}}^+ (\tg) \subset sp_{{\cal F}}^+ (\tf)\subset
\Lambda $. To complete the proof of the theorem we need to show
that $\lim_{\alpha \downarrow 0} \alpha R(\alpha +i\xi ,\tD )\tg
=\tilde{0}$ for all $\xi\in\Lambda$. In fact, for each $Re\lambda
\not= 0$, we have
\begin{eqnarray*}
\lim_{\alpha \downarrow 0} \alpha R(\alpha +i\xi ,\tD )R(\lambda
,{\tD}) \tf  &=& \lim_{\alpha \downarrow 0} R(\lambda ,
 {\tD})\alpha R(\alpha
+i\xi ,\tD )\tf\\
&=& R(\lambda ,
 {\tD}) \lim_{\alpha \downarrow 0} \alpha R(\alpha +i\xi ,\tD )\tf
 =0 .
\end{eqnarray*}
By the above argument used to show (\ref{erg}), this shows that
$\lim_{\alpha \downarrow 0} \alpha R(\alpha +i\xi ,\tD )\tg =0$.
\end{proof}
\begin{theorem}\label{the loomis half line}
Let ${\cal F}$ be a function space of $BC(\R^+,\X)$ that satisfies
Condition $F^+$, and let $f$ be in $BC(\R^+,\X)$ such that
$sp^+_{\cal F}(f)$ is countable. Moreover, assume that
\begin{equation}
\lim_{\alpha \downarrow 0} \alpha R(\alpha +i\xi ,\tD )\tf =0
\end{equation}
for all $\xi \in sp^+_{\cal F}(f)$. Then, $f\in {\cal F}$.
\end{theorem}
\begin{proof}
Set $\Lambda := sp^+_{\cal F}(f)$. Consider the function space $
{\Lambda}^e_{\cal F}(\X)$ and the operator $\tD_\Lambda $ on it.
We are going to prove that the function space ${\Lambda}^e_{\cal
F}(\X)$ is trivial. In fact, let us assume to the contrary that it
is not trivial. Then, since $sp^+_{\cal F}(f)$ is countable, by
Proposition \ref{the spe dif ope hl}, the spectrum $\sigma
(\tD_\Lambda )$ is countable. Therefore, there is an isolated
point. By Lemma \ref{lem 4} and Theorem \ref{the gel}, this isolated
point of spectrum of $\sigma (\tD_\Lambda )$ must be an
eigenvalue. And hence, there exists a nonzero vector $\tilde{g}\in
{\Lambda}^e_{\cal F}(\X)$ such that $(\tD -i\xi )\tilde{g}=0$. For
$\alpha >0$, since $(\tD -i\xi )\tilde{g}=0$ by (\ref{3.13}) we
have $
 \alpha   R(\alpha +i\xi ,\tD ) \tilde{g}   =
 \tilde{g} .
$ As $\tg\in {\Lambda}^e_{\F}(\X)$,
\begin{eqnarray*}
0 = \lim_{\alpha \to 0^+} \alpha   R(\alpha +i\xi ,\tD ) \tilde{g}
= \lim_{\alpha \to 0^+} \tilde{g} =  \tilde{g} .
\end{eqnarray*}
This contradiction shows that ${\Lambda}^e_{\cal F}(\X)$ must be
trivial. Therefore, $f\in {\cal F}$.
\end{proof}
\begin{remark}
For $f\in \F \subset BUC(\R^+,\X)$, there is a relation between
$sp_\F ^+(f)$ and the set $Sp^+(f)$ of all singularities of the
Laplace transform of $f$, that is, the set of all reals $\xi$ such
that the Laplace transform $\hat{f}(\lambda )$ of $f$ has no
analytic extension to any neighborhood of $i\xi$. In fact, if
$f\in \F \subset BUC(\R^+,\X)$, then
\begin{equation}
sp_\F ^+(f) \subset Sp^+(f) .
\end{equation}
This is a consequence of \cite[Theorem 5.3.4, p. 171]{nee}. So,
Theorem \ref{the loomis half line} extends  \cite[Theorem
2.3]{arebat2} and \cite[Theorem 4.1]{batneerab}.
\end{remark}
The following corollaries follow immediately from Theorem \ref{the
loomis half line}:
\begin{corollary}
Let $f\in BC(\R^+,\X)$ such that $sp^+_{AAP(\X)}(f)$ is countable.
Moreover, assume that
\begin{equation}
\lim_{\alpha \downarrow 0} \alpha R(\alpha +i\xi ,\tD )\tf =0
\end{equation}
for all $\xi \in sp^+_{AAP(\X)}(f)$. Then, $f$ is asymptotically
almost periodic.
\end{corollary}
\begin{corollary}
Let $f\in BC(\R^+,\X)$ such that $sp^+_{AA(\X)}(f)$ is countable.
Moreover, assume that
\begin{equation}
\lim_{\alpha \downarrow 0} \alpha R(\alpha +i\xi ,\tD )\tf =0
\end{equation}
for all $\xi \in sp^+_{AA(\X)}(f)$. Then, $f$ is asymptotically
almost automorphic.
\end{corollary}

\section{Applications to the Asymptotic Behavior of Solutions of Evolution Equations}
\subsection{Equations on the whole line}
Consider evolution equations of the form
\begin{equation}\label{app 1}
\frac{du(t)}{dt} = Au(t) +f(t), \quad t\in \R , \ u(t)\in \X ,
\end{equation}
where $A$ is a closed linear operator on a Banach space $\X$, $f$
is an $\X$-bounded and continuous function on $\R$. Throughout
this section we always assume that $A$ is a closed linear
operator.
\begin{definition}
A function $u\in BC(\R,\X)$ is said to be a mild solution of
(\ref{app 1}) if for every $t\in \R$, $\int^t_0 u(s)ds \in D(A)$,
and
\begin{equation}\label{mild sol}
u(t) -u(0) = A\int^t_0 u(s)ds + \int^t_0 f(s)ds , \quad \
\mbox{for all} \ t\in \R .
\end{equation}
\end{definition}
The following lemma and its proof have been known in the uniform
continuity setting (see e.g.
\cite{arebathieneu,delvu,levzhi,bas}). For the reader's
convenience we re-state its version for non-uniform continuous
mild solutions with a standard proof.
\begin{lemma}\label{lem sp cal F of u}
Let ${\cal F}\subset BC(\R,\X)$ be a function space that satisfies
Condition F, and $f\in BC(\R,\X)$, and let $u\in BC(\R,\X)$ be a
mild solution of (\ref{app 1}) on $\R$. Then,
\begin{equation}
sp_{\cal F}u \subset \sigma_i(A) \cup sp_\F (f),
\end{equation}
where $\sigma_i(A):= \{ \xi \in \R :\ i\xi \in \sigma (A)\}$.
\end{lemma}
\begin{proof}
For every $Re\lambda >0$, and $s \in \R$, we have
$$
\int^\infty _0 e^{-\lambda t} \left(\int^{t+s}_0 u(\xi )d\xi
\right) dt = \frac{1}{\lambda}\left( \int^\infty_0 e^{-\lambda
t}u(t+s )dt +\int^s_0 u(\xi )d\xi \right).
$$
Applying this to (\ref{mild sol}), for every $Re\lambda >0$, and
$s \in \R$, since
$$
u(s) = u(0) + \int^s_0 u(\xi )d\xi +\int^s_0 f(\xi) d\xi , \quad
s\in \R
$$
we have
\begin{eqnarray*}
\int^\infty _0 e^{-\lambda t} u(t+s) dt &=& \int^\infty _0
e^{-\lambda t} dt u(0) + A \int^\infty _0 e^{-\lambda t}
\left(\int^{t+s}_0 u(\xi )d\xi \right) dt \\
&&\hspace{3cm} + \int^\infty _0 e^{-\lambda t} \left(\int^{t+s}_0
f(\xi )d\xi \right) dt \\
&=& \frac{1}{\lambda} u(0) + \frac{1}{\lambda} A \int^\infty _0
e^{-\lambda t} u(t+s) dt + \frac{1}{\lambda} A \int^s_0 u(\xi)d\xi
\\
&& \hspace{3cm} \frac{1}{\lambda}  \int^\infty _0 e^{-\lambda t}
f(t+s) dt + \frac{1}{\lambda}  \int^s_0 f(\xi)d\xi \\
&=&  \frac{1}{\lambda} A \int^\infty _0 e^{-\lambda t} u(t+s) dt +
\frac{1}{\lambda}  \int^\infty _0 e^{-\lambda t} f(t+s) dt +
\frac{1}{\lambda}u(s) .
\end{eqnarray*}
Therefore, for $Re\lambda >0$, by (\ref{re-1}),
\begin{equation}\label{spe of mild sol}
R(\lambda , \cD) u= \frac{1}{\lambda} {\cal A}R(\lambda ,\cD) u
+\frac{1}{\lambda}  R(\lambda ,\cD) f +\frac{1}{\lambda}u ,
\end{equation}
where $\cal A$ denotes the operator of multiplication by $A$ on
$BC(\R,\X)$. Similarly, we can show that (\ref{spe of mild sol})
holds also for $Re\lambda <0$. Therefore, for $Re\lambda \not= 0$,
\begin{equation}\label{spe of mild sol 2}
(\lambda -{\cal A})R(\lambda , \cD) u=   R(\lambda ,\cD) f +u ,
\end{equation}
By the axioms defining Condition F, and the assumption, we arrive
at
\begin{eqnarray}
(\lambda -\tilde{{\cal A}})R(\lambda , \tD) \tilde{u} &=&
R(\lambda ,\tD) \tf +\tilde{u}  .
\end{eqnarray}
If $\xi_0\in \R\backslash \sigma_i(A)$ and $\xi_0\not\in sp_\F
(f)$, for $\lambda$ in a small neighborhood of $i\xi_0$, and
$Re\lambda\not= 0$,
\begin{eqnarray}
R(\lambda , \tD) \tilde{u} &=& R(\lambda ,\tilde{{\cal
A}})R(\lambda , \tD)\tilde{f} + R(\lambda ,\tilde{{\cal
A}})\tilde{u} .
\end{eqnarray}
Therefore, $R(\lambda , \tD) \tilde{u}$ has an analytic extension
to a neighborhood of $i\xi_0$, so $\xi_0\not\in sp_{\cal F}(u)$.
This proves the lemma.
\end{proof}
The following corollary is an immediate consequence of the above
lemma and Theorem \ref{the loomis}.
\begin{corollary}\label{cor 4.3}
Let ${\cal F}\subset BC(\R,\X)$ be a function space that satisfies
Condition F and contains $f$, and let $u\in BC(\R,\X)$ be a mild
solution of (\ref{app 1}) on $\R$. Moreover, assume that
$\sigma_i(A)$ be countable. Then, $u$ is in ${\cal F}$.
\end{corollary}
\begin{remark}
When ${\cal F}\subset BUC(\R,\X)$ and $u\in BUC(\R,\X)$, the above
corollary is known in \cite{arebat,bas,bask,bask2,levzhi,ruevu}
that extends Loomis Theorem for the scalar functions to vector
valued ones. In these works the assumption on the uniform
continuity is essential to make use of the techniques based on the
spectral properties of $C_0$-groups.
\end{remark}
The following are standard corollaries of Corollary \ref{cor 4.3}.
\begin{corollary}
Let $f\in AP(\X)$, $\X$ not contain any subspace isomorphic to
$c_0$, and let $u\in BC(\R,\X)$ be a mild solution on $\R$ of
(\ref{app 1}) for which $\sigma (A)\cap i\R $ is countable. Then,
$u$ is almost periodic.
\end{corollary}
\begin{corollary}
Let $f\in AA(\X)$,  $\X$ not contain any subspace isomorphic to
$c_0$, and let $u\in BC(\R,\X)$ be a mild solution on $\R$ of
(\ref{app 1}) for which $\sigma (A)\cap i\R $ is countable. Then,
$u$ is almost automorphic.
\end{corollary}

\subsection{Equations on the half - line}
In this subsection we consider linear evolution equations on the
half line, that is,
\begin{equation}\label{eq hl}
\frac{du(t)}{dt} = Au(t) +f(t), \quad t\in \R^+ , \ u(t)\in \X ,
\end{equation}
where $f\in BC(\R^+,\X)$, $A$ is a closed linear operator on $\X$.

\begin{lemma}\label{lem spe mil hl}
Let ${\cal F}\subset BC(\R^+,\X)$ be a function space that
satisfies Condition $F^+$ and contains $f$, and let $u\in
BC(\R^+,\X)$ be a mild solution of (\ref{eq hl}) on $\R^+$. Then,
\begin{equation}\label{est of spe of mil hl}
sp^+_{\cal F}u \subset \sigma_i(A).
\end{equation}
Moreover, for $Re\lambda \not= 0$,
\begin{eqnarray}\label{4.12}
R(\lambda , \tD) \tilde{u} &=& R(\lambda ,\tilde{{\cal
A}})\tilde{u} .
\end{eqnarray}
\end{lemma}
\begin{proof}
Let $Re\lambda >0$. In the same way as in the proof of Lemma
\ref{lem sp cal F of u}, by (\ref{resolvent of tilde D})
\begin{equation}\label{est of sp cal F hl}
(\lambda -\tilde{{\cal A}})R(\lambda , \tD) \tilde{u}  =
R(\lambda ,\tD) \tf +\tilde{u} = \tilde{u} , \quad Re\lambda >0.
\end{equation}
We now show that (\ref{est of sp cal F hl}) holds for $Re\lambda
<0$ as well. For each $Re\lambda <0$, using the
integration-by-parts formula we have
\begin{eqnarray*}
\int^t_0 e^{\lambda (t-s)}\left(  \int^s_0 u(\xi )d\xi \right)  ds
&=& -\frac{e^{\lambda(t-s)}}{\lambda} \int^s_0 u(\xi )d\xi
\mid^t_0 + \frac{1}{\lambda} \int^t_0 e^{\lambda (t-s)}u(\xi)d\xi
\\
&=& -\frac{1}{\lambda} \int^t_0 u(\xi )d\xi + \frac{1}{\lambda}
\int^t_0 e^{\lambda (t-s)}u(\xi)d\xi .
\end{eqnarray*}
Applying this to (\ref{mild sol}),  we arrive at
\begin{eqnarray*}
 \int^t_0 e^{\lambda(t-\xi)}u(\xi )d\xi &=&  \int^t_0 e^{\lambda(t-\xi)}u(0 )d\xi
+ A  \int^t_0 e^{\lambda (t-s)}\left(  \int^s_0 u(\xi )d\xi \right)  ds   \\
&&\hspace{2cm} + \int^t_0 e^{\lambda (t-s)}\left(  \int^s_0 f(\xi )d\xi \right)  ds \\
\int^t_0 e^{\lambda(t-\xi)}u(\xi )d\xi  &=&
 \frac{e^{\lambda
t}u(0)}{\lambda} -\frac{u(0)}{\lambda}+A\left( -\frac{1}{\lambda}
\int^t_0 u(\xi )d\xi +
\frac{1}{\lambda} \int^t_0 e^{\lambda (t-s)}u(\xi)d\xi \right) \\
&&\hspace{2cm}  -\frac{1}{\lambda} \int^t_0 f(\xi )d\xi +
\frac{1}{\lambda} \int^t_0 e^{\lambda (t-s)}f(\xi)d\xi .
\end{eqnarray*}
Note that for each $Re\lambda <0$,
\begin{eqnarray*}
-\frac{u(0)}{\lambda} -\frac{1}{\lambda} A \int^t_0 u(\xi )d\xi
-\frac{1}{\lambda} \int^t_0 f(\xi )d\xi = -\frac{u(t)}{\lambda},
\quad \mbox{for all } s \in\R^+ .
\end{eqnarray*}

 Therefore,
\begin{eqnarray*}
 \int^t_0 e^{\lambda(t-\xi)}u(\xi )d\xi &=& \frac{1}{\lambda} A
 \int^t_0 e^{\lambda(t-\xi)}u(\xi )d\xi + \frac{1}{\lambda} \int^t_0 e^{\lambda(t-\xi)}f(\xi
 )d\xi\\
 &&\hspace{2cm}-\frac{1}{\lambda}u(t) + \frac{e^{\lambda
 t}u(0)}{\lambda}.
\end{eqnarray*}
Note that
$$
 \R^+ \ni t \mapsto \frac{1}{\lambda} e^{\lambda t}u(0) \in \X \
\ \mbox{belongs to} \ {\cal F}.
$$
Therefore, by (\ref{resolvent of tilde D}), for $Re\lambda <0$ we
have
\begin{equation}\label{spe of mild sol hl}
-R(\lambda , \tD) \tilde{u}= -\frac{1}{\lambda} {\cal A}R(\lambda
,\tD) \tilde{u} -\frac{1}{\lambda}  R(\lambda ,\tD)
\tf-\frac{1}{\lambda}\tilde{u} .
\end{equation}
So, (\ref{est of sp cal F hl}) holds for $Re\lambda <0$ as well.
Next, if $\xi_0\in \R\backslash \sigma_i(A)$, then for $\lambda$
in a sufficiently small neighborhood of $i\xi_0$,
\begin{eqnarray}
R(\lambda , \tD) \tilde{u} &=& R(\lambda ,\tilde{{\cal
A}})\tilde{u} .
\end{eqnarray}
Therefore, $R(\lambda , \tD) \tilde{u}$ has an analytic extension
to a neighborhood of $i\xi_0$, so $\xi_0\not\in sp_{\cal F}(u)$.
This proves the lemma.
\end{proof}

\begin{remark}
For $f,u\in BUC(\R^+,\X)$, the estimate (\ref{est of spe of mil
hl}) has been made in \cite{arebat2} by a different method that
seems to be unapplicable to the case of non-uniformly continuous
functions $f$ and $u$. In our approach to the spectrum, the
resolvent $R(\lambda , \tD )\tf$ for $Re\lambda <0$ is explicitly
found. So, the above lemma can be proved much easier than in the
previous works.
\end{remark}
Let $f\in BC(\R^+,\X)$ be uniformly ergodic at $i\xi$ for some
$\xi\in\R$, that is, the following limit exists
$$
 \lim_{\alpha \downarrow 0} \alpha R(\alpha +i\xi , \cD)f = \bar f
 \in BC(\R^+,\X) .
$$
 Since $g:= R(\alpha +i\xi , \cD)f$
satisfies the equation $(\alpha +i\xi )g(t) - g'(t) = f(t)$, by
the Variation-of-Constants Formula,
$$
g(t) = e^{(\alpha +i\xi )t}g(0) - \int^t_0 e^{(\alpha +i\xi
)(t-s)}f(s)ds, \quad \alpha >0, t\in \R^+ .
$$
Therefore, for every fixed $t\in\R^+$,
\begin{eqnarray*}
\bar f(t) &=& \lim_{\alpha \downarrow 0}\alpha g(t) = \lim_{\alpha
\downarrow 0}\alpha e^{(\alpha +i\xi )t}g(0) - \lim_{\alpha
\downarrow 0}\alpha \int^t_0 e^{(\alpha +i\xi
)(t-s)}f(s)ds\\
&=& \lim_{\alpha \downarrow 0}\alpha e^{(\alpha +i\xi )t}g(0)\\
&=& e^{i\xi t} \bar f(0) .
\end{eqnarray*}
This shows that if a function $f\in BC(\R^+,\X)$ is uniformly
ergodic at $i\xi$ for some $\xi\in\R$, and
$$
 \lim_{\alpha \downarrow 0} \alpha R(\alpha +i\xi , \cD)f = \bar f
 \in BC(\R^+,\X) ,
$$
then for each $t\in\R^+$, $\bar f(t) = e^{i\xi t}a$ for some fixed
$a\in\X$, so $\bar f\in AAP(\R^+,\X)$.

 The following corollaries are obvious.
\begin{corollary}
Let $\sigma_i(A)$ be countable, and let $f$ be asymptotically
almost periodic. Then every bounded mild solution $u$ on $\R^+$ of
(\ref{eq hl}) is asymptotically almost periodic provided $u$ is
uniformly ergodic at $i\xi$ for each $\xi\in\sigma _i(A)$.
\end{corollary}
\begin{corollary}
Let $\sigma_i(A)$ be countable, and let $f$ be asymptotically
almost automorphic. Then every bounded mild solution $u$ of
(\ref{eq hl}) on $[0,\infty )$ is asymptotically almost
automorphic provided $u$ is uniformly ergodic at $i\xi$ for each
$\xi\in\sigma _i(A)$.
\end{corollary}

\medskip
Let us consider the homogeneous equation
\begin{equation}\label{eq hl homo}
\dot u(t) =Au(t),\quad u(t)\in \X , \ t\in \R^+ ,
\end{equation}
where $A$ is a closed linear operator on $\X$. A mild solution $u$
on $\R^+$ of (\ref{eq hl homo}) is {\it asymptotically stable} if
$\lim_{t\to\infty}u(t)=0$.

\begin{theorem}\label{the minh}
Let $u\in BC(\R^+,\X)$ be a mild solution of (\ref{eq hl homo}),
and let $A$ satisfy the following conditions:
\begin{enumerate}
\item $ \sigma  (A) \cap i\R $ is countable; \item $ \lim_{\alpha
\downarrow 0}   \alpha R(\alpha +i\xi  , A)u(t)   =0 $ uniformly
in $t\in\R^+$,  for all $ i\xi \in \sigma (A)\cap i\R $.
\end{enumerate}
Then, the solution $u$ is asymptotically stable.
\end{theorem}
\begin{proof}
Let $u$ be a bounded mild solution on $\R^+$ of (\ref{eq hl
homo}), and let $\F := C_0(\R^+,\X)$. Then by Lemma \ref{lem spe
mil hl}, $sp^+_{C_0(\R^+,\X)}(f)$ is countable. Moreover, by
(\ref{4.12}), for $\alpha > 0$, we have
\begin{equation}
\| R(\alpha +i\xi , \tD) \tilde{u} \|  = \|  R(\alpha +i\xi
,\tilde{{\cal A}})\tilde{u}\| .
\end{equation}
Therefore, for all $\xi \in \sigma_i(A)$.
\begin{equation}
0 \le \lim_{\alpha \downarrow 0} \| \alpha R(\alpha +i\xi , \tD)
\tilde{u} \|  = \lim_{\alpha \downarrow 0} \| \alpha R(\alpha
+i\xi , \tilde{{\cal A}}) \tilde{u} \| = 0.
\end{equation}
Applying Theorem \ref{the loomis half line}, we end up with $u\in
C_0(\R^+,\X)$, proving the theorem.
\end{proof}
The following corollary is an immediate consequence of Theorem
\ref{the minh}.

\begin{corollary}\label{cor minh}
Let $A$ satisfy the following conditions:
\begin{enumerate}
\item $ \sigma  (A) \cap i\R $ is countable; \item $ \lim_{\alpha
\downarrow 0}   \alpha R(\alpha +i\xi  , A)   =0 $ for all $ i\xi
\in \sigma (A)\cap i\R $.
\end{enumerate}
Then, every bounded mild solution on $\R^+$ of (\ref{eq hl homo})
is asymptotically stable.
\end{corollary}

\begin{remark}
If $ \sigma  (A) \cap i\R =\emptyset$, the condition (ii) in the
above theorem follows immediately from the condition (i).
\end{remark}
If $A$ is the infinitesimal generator of a bounded
$C_0$-semigroup, then, as a consequence of Theorem \ref{the minh}
we obtain the following well-known Arendt-Batty-Lyubich-Vu
Theorem.

\begin{corollary}\label{the ablv}
{\rm (The Arendt-Batty-Lyubich-Vu Theorem)} Let $A$ generate a
bounded $C_0$-semigroup on a Banach space $\X$, and let it satisfy
the following conditions:
\begin{enumerate}
\item $ \sigma  (A) \cap i\R $ is countable; \item $ \sigma_p(A^*)
\cap i\R =\emptyset  $.
\end{enumerate}
Then, every mild solution on $\R^+$ of (\ref{eq hl homo}) is
asymptotically stable.
\end{corollary}
\begin{proof}
Since $A$ generates a $C_0$-semigroup $(T(t))_{t\ge 0}$, each mild
solution $u$ is of the form $u(t)=T(t)x$, for all $t\in\R^+$ and
some $x\in \X$, so $u\in BC(\R^+,\X)$. As is well known (see e.g.
\cite[Sect. 5.5]{arebathieneu}, \cite[Proposition
3.2]{batneerab2}) the conditions (ii) yields the following
$$
\lim_{\alpha \downarrow 0}   \alpha R(\alpha +i\xi  , A) x
 = 0, \ \mbox{for all} \ \ x,\in \X , i\xi \in \sigma (A)\cap i\R
 .
$$
Therefore, for all  $ \ \alpha >0 , i\xi \in \sigma (A)\cap i\R $,
\begin{eqnarray*}
\lim_{\alpha\downarrow 0}  \sup_{t\in\R^+}\| \alpha R(\alpha +i\xi
,A)u(t)\|
& =& \lim_{\alpha\downarrow 0} \sup_{t\in\R^+}\| \alpha R(\alpha +i\xi ,A)T(t)x \| \\
&=&\lim_{\alpha\downarrow 0}  \sup_{t\in\R^+} \| \alpha T(t)
\int^\infty _0 e^{-(\alpha +i\xi )s}
T(s)x dt\|  ,\\
&\le& \lim_{\alpha\downarrow 0} \sup_{t\in\R^+}  \|  T(t)\| \cdot
\alpha\|\int^\infty _0 e^{-(\alpha +i\xi )s} T(s)x dt\| \\
&=& 0 .
\end{eqnarray*}
So, by Theorem \ref{the minh}, $u$ is asymptotically stable.
\end{proof}
\begin{remark}
The Arendt-Batty-Lyubich-Vu Theorem was proved independently by
Arendt and Batty in \cite{arebat3}, and Lyubich and Vu in
\cite{lyuvu}. Earlier in \cite{sklshi} Sklyar\footnote{The author
thanks G. M. Sklyar for sending him a copy of the original paper
\cite{sklshi}} and Shirman proved a similar result for bounded $A$
using a method based on the concept of "isometric limit
semigroups" which can be extended to the case where $A$ is the
generator of a $C_0$-semigroup. There are many extensions of the
Arendt-Batty-Lyubich-Vu Theorem (see e.g.
\cite{batneerab,batneerab2,arebat2,delvu}). Note that in all these
extensions the assumption on the uniform continuity of mild
solutions is essential due to the techniques using the
 theory of $C_0$-semigroups. If $A$ generates a $C_0$-semigroup,
 the uniform continuity of mild solutions on $\R^+$ follows from the condition of (ii) in the
 above corollary.
\end{remark}

\end{document}